\documentclass[12pt]{article}

\usepackage[T2A]{fontenc}
\usepackage[utf8]{inputenc}

\usepackage{amsfonts}
\usepackage{amssymb}
\usepackage{amsmath}
\usepackage{amsthm}

\def \le {\leqslant}
\def \ge {\geqslant}

\begin{document}

\begin{Large}
\centerline{\bf
On Diophantine approximations with positive integers:
}
\centerline{\bf
a remark to W.M.Schmidt's theorem
  }
 
 \vskip+1.5cm \centerline{by {\bf Nicky Moshchevitin}\footnote{ Research is supported by
the grant RFBR No. 09-01-00371-a}} \vskip+1.5cm
\end{Large}
\vskip+0.7cm
\centerline{\bf Abstract.}
\begin{small}
 We prove a generalization of W.M. Schmidt's theorem related to the Diophantine approximations for a linear form of the type
$\alpha_1x_1+\alpha_2x_2 +y$ with {\it positive} integers $x_1,x_2$.
\end{small}
\vskip+2.0cm

\section{Introduction}
Let $||\xi ||$ denotes the distance from real $\xi$ to the nearest integer.
Let
$\tau =\frac{1+\sqrt{5}}{2}
.$
In \cite{SCH}
W.M. Schmidt proved the following result.
\vskip+0.5cm
{\bf Theorem 1.} (W.M.Schmidt)\,\,{\it
Let real numbers $\alpha_1,\alpha_2$
be linearly independent over $\mathbb{Z}$
together with 1. Then there exists a sequence of integer two-dimensional vectors
 $(x_1(i), x_2(i))$
 such that
}

1.\,\, $x_1(i), x_2(i) > 0$;

2.\,\, $||\alpha_1x_1(i)+\alpha_2x_2(i) ||\cdot (\max \{x_1(i),x_2(i)\})^\tau \to 0$ as $ i\to +\infty$.
\vskip+0.5cm

A famous conjecture that the exponent $\tau$ here may be  replaced by $2-\varepsilon$ with arbitrary positive $\varepsilon$
(see \cite{SCH,SCH1}) is still unsolved.
We would like to mention that there are various generalizations of 
W.M. Schmidt's theorem  by P.Thurnheer \cite{T1,T2} Y. Bugeaud and S. Kristensen \cite{BK} and some other mathematicians.
 
For a real $\gamma \ge 2$ we define a function
$$
g(\gamma ) = \tau +\frac{2\tau - 2}{\tau^2\gamma -2}.
$$
One can see that $g(\gamma )$ is a strictly decresaing function and
$$
g(2) =2,\,\,\,\,\,\, \lim_{\gamma\to+\infty} g(\gamma ) = \tau.
$$
For positive $\Gamma$ define
$$
C(\Gamma ) = 2^{18} \Gamma^{\frac{\tau -\tau^2}{\tau^2\gamma-2}}.
$$

In this paper we prove  the following statement. 
\vskip+0.5cm
{\bf Theorem 2.} \,\,{\it
Suppose that real numbers $\alpha_1,\alpha_2$
satisfy the following Diophantine condition. 
For some $\Gamma\in (0,1)$ and $\gamma \ge 2$ the inequality
\begin{equation}
||\alpha_1m_1+\alpha_2m_2||\ge \frac{\Gamma}{(\max \{|m_1|,|m_2|\})^\gamma}
\label{bad}
\end{equation}
holds for all integer vectors $(m_1,m_2)\in \mathbb{Z}^2\setminus \{(0,0)\}$.
Then
 there exists an infinite sequence of integer two-dimensional vectors
 $(x_1(i), x_2(i))$
 such that}

1.\,\, $x_1(i), x_2(i) > 0$;

2.\,\, $||\alpha_1x_1(i)+\alpha_2x_2(i) ||\cdot (\max \{x_1(i),x_2(i)\})^{g(\gamma )} \le C(\Gamma )$ for all  $ i$.
\vskip+0.5cm
Of course  the constant $2^{18}$ in the definition of $C(\Gamma )$ may be reduced.

\section{The best approximations}

 Suppose that $1,\alpha_1,\dots,\alpha_r$ are
linearly independent over $\mathbb{Z}$. For an integer point ${\bf m}=(m_0,m_1,m_2)\in\mathbb{Z}^3$ we
define
$$
\zeta({\bf m})=m_0+m_1\alpha_1+m_2\alpha_2.
$$
A point ${\bf m}=(m_0,m_1,m_2)\in\mathbb{Z}^{3}\setminus\{(0,0,0)\}$ is defined to be \emph{a best
approximation (in the sense of linear form)}  if
$$
  \zeta({\bf m})=\min_{\bf n}\|\zeta({\bf n})\|,
$$
where  the minimum is taken over all the
integer vectors ${\bf n}=(n_0,n_1,n_2)\in\mathbb{Z}^3$ such that
$$
0<\max_{1\le j\le 2 }|n_j|\le \max_{1\le  j\le  2} |m_j|.
$$
All the best approximations form a sequence of points ${\bf m}_\nu=(m_{0,\nu},m_{1,\nu},m_{2,\nu})$
with increasing $\max_{1\le  j\le  2}|m_{j,\nu}|$.  

Let us denote
$$
  \zeta_\nu=\zeta({\bf m}_\nu),\quad M_\nu=\max_{1\le  j\le  2}|m_{j,\nu}|.
$$
Then
$$
\zeta_1>\zeta_2>\cdots>\zeta_\nu>\zeta_{\nu+1}>\cdots
$$
and
$$
M_1<M_2<\cdots<M_\nu<M_{\nu +1}<\cdots. 
$$
It follows from the Minkowski convex body theorem that 
\begin{equation}
\zeta_\nu M_{\nu+1}^2\le 1.
\label{mink}
\end{equation} 
To  prove the inequality (\ref{mink}) one should consider the parallelepiped $\Omega_\nu \in \mathbb{R}^3$
which consisit of all points  $(y,x_1,x_2)\in \mathbb{R}^3$ satisfying the inequalities
$$
\begin{cases}
|y+\alpha_1 x_1+\alpha_2x_2|<\zeta_{\nu},
\cr
\max\{ |x_1|,|x_2|\} < M_{\nu+1}
\end{cases}
.$$
Then there is no non-zero integer points in $\Omega_\nu$ and (\ref{mink}) follows.

\section{A statement about consecutive best approximations}
Now we formulate a rather technical result.
\vskip+0.5cm

{\bf Theorem 3.}\,\,{\it
Let   
\begin{equation}
\left|
\begin{array}{ccc}
m_{0,\nu-1} &m_{1,\nu -1} & m_{2,\nu-1}\cr
m_{0,\nu} &m_{1,\nu } & m_{2,\nu}\cr
m_{0,\nu+1} &m_{1,\nu +1} & m_{2,\nu+1}
\end{array}
\right|\neq 0.
\label{neq}
\end{equation}
Then at least one of two satements below is valid. 
\vskip+0.3cm
{\bf (i)}
\,\,
There exists an integer point
 $(x_1^0,x_2^0)$ 
such that}

1.\,\, $x_1^0, x_2^0 >0$;

2.\,\, $ M_{\nu+2}\le \max \{x_1^0,x_2^0 \} \le  4M_{\nu+2}$;

3.\,\, $||\alpha_1x_1^0+\alpha_2 x_2^0 ||\le 16 (\max \{x_1^0,x_2^0\})^{-2}.$

\vskip+0.3cm
{\bf (ii)}\,\,{\it  There exists an integer point $(x_1^0,x_2^0)$  such that
}

1.\,\, $x_1^0, x_2^0 >0$;

2.\,\, $\max \{x_1^0,x_2^0 \} \le  240M_{\nu+1}^{\tau}M_\nu^{-\frac{1}{\tau}}$;

3.\,\, $||\alpha_1x_1^0+\alpha_2 x_2^0 ||\le24^\tau M_\nu^{\frac{1-\tau}{\tau}} (\max\{ x_1^0,x_2^0\})^{-\tau}.$

\vskip+0.5cm
It is a well-known fact (see for example \cite{DS,M})  that
there exists infinitely many $\nu$ such that (\ref{neq}) holds, provided that the numbers
  $1, \alpha_1,\alpha_2$ are linearly independent over $\mathbb{Z}$.
So Theorem 1 follows from Theorem 3 as $M_\nu \to +\infty$,  $\nu\to +\infty$.

Here we would like to give few comments.
Theorem 3 may be treated as a "local" statement which provides the existence of
a small  value of the linear form $||\alpha_1x_1+\alpha_2x_2||$
relatively "close"  to the best approximations satisfying
  (\ref{neq}). 
We shall give the proof of Theorem 3 in next two sections. The proof follows the original construction due to W.M.Schmidt
\cite{SCH}, however it includes few modifications.
\vskip+0.5cm
Now we show that Theorem 3 implies Theorem 2.

Suppose that  the statement {\bf (i)} holds for infinitely many $\nu$. Then 
as
$C(\Gamma ) \ge 16,\,\, g (\gamma )\le 2$
we see that Theorem 2 follows from  Theorem 3
obviously.

So we may assume that the statement {\bf (ii)} holds for infinitely many $\nu$.
From the condition (\ref{bad}) of Theorem 2 and from the inequality (\ref{mink}) 
applied to the vector $(m_1,m_2) =(m_{1,\nu},m_{2,\nu})$
we deduce that
$$
\Gamma M_\nu^{-\gamma}\le M_{\nu+1}^{-2}.
$$
The last inequality together with the statement 2 of {\bf (ii)} gives
$$
240^{-\frac{2\tau}{\tau^2\gamma -2}}\times
\Gamma^{\frac{\tau^2}{\tau^2\gamma -2}}\times
\left(\max \{x_1^0,x_2^0 \} \right)^{\frac{2\tau}{\tau^2\gamma -2}} \le M_\nu.
$$
Now we substitute the last inequality into the statement 3 of {\bf (ii)} and obtain 
$$
||\alpha_1x_1(i)+\alpha_2x_2(i) ||\le \frac{C(\Gamma )}{ (\max \{x_1(i),x_2(i)\})^{g(\gamma )}}
.
$$
Theorem 2 is proved.

\section{Lemmata}

Put  
$$R_\nu = 2(M_{\nu+1}\zeta_{\nu})^{-1}.$$
 From (\ref{mink}) it follows that
$
R_\nu > M_{\nu+1}.
$
\vskip+0.5cm

{\bf Lemma 1.}\,\,{\it 
Let numbers $1,\alpha_1,\alpha_2$, be linearly independent over  $\mathbb{Z}$.
Then there exists an integer point 
${\bf x}^0= (x_1^0,x_2^0)$
such that
}

1.\,\, $x_1^0, x_2^0 >0$;

2.\,\, $\max \{x_1^0,x_2^0 \} \le R_\nu$;

3.\,\, $||\alpha_1x_1^0+\alpha_2 x_2^0 ||<\zeta_{\nu}.$

Proof.

 Consider the papallelepiped
$\Omega_\nu^{1}$ defined by the system of inequalities
$$
\begin{cases}
|\alpha_1 x_1+\alpha_2x_2+y|
\le\zeta_{\nu},
\cr
|x_1-x_2|\le M_{\nu+1},
\cr
|x_1+x_2 |\le R_\nu
.
\end{cases}
$$
 As $M_{\nu+1} R_\mu \zeta_{\nu} = 2$, the measure of $\Omega_\nu^{1}$ is equal to
   8.
Hence by the Minkowski convex body theorem there exists a non-zero integer point 
${\bf z}^0 = (y^0,x_1^0,x_2^0)\in \mathbb{Z}^3\cap \Omega_\nu^{1}$.
As it was mentioned in Section 2 parallelepiped
 $\Omega_\nu$ contains no non-zero integer points.
So
$ 
{\bf z}\in \Omega_2\setminus \Omega_1$.
We see that for the integers $x_1^0,x_2^0$ the statements 1 -3 of Lemma 1 are true
(the strict inequalities in 1 and 3 follow from the linear independence of $1,\alpha_1,\alpha_2$).

Lemma is proved.
 
\vskip+0.5cm
{\bf Remark.}\,
As the inequality in the statement 3 of Lemma 1 is a strict one we deduce that
 $ \max\{x_1^0,x_2^0\} \ge M_{\nu+1}$.
\vskip+0.5cm

{\bf Corollary 1.}\,\,{\it
Let the following inequality be valid:
\begin{equation}
\zeta_{\nu } \ge (8 M_{\nu+1}^{2})^{-1}
.
\label{oche}
\end{equation}
 Then there exists
an integer point ${\bf x}^0= (x_1^0,x_2^0)$ such that
}

1.\,\, $x_1^0, x_2^0 >0$;

2.\,\, $M_{\nu+1}\le \max \{x_1^0,x_2^0 \} \le  4M_{\nu+1}$;

3.\,\, $||\alpha_1x_1^0+\alpha_2 x_2^0 ||\le 16 (\max \{x_1^0,x_2^0\})^{-2}.$

\vskip+0.5cm

Proof.

Apply Lemma 1. The numbers
 $x_j^0$  from Lemma 1 are positive.
Inequality (\ref{oche}) and the remark after Lemma 1 lead to the statement 2 of Corollary 1.
Now we apply the statement  3 of Lemma 1, the inequality  (\ref{mink}) and the statement 2 of Corollary 1 to see that
 $$
||\alpha_1x_1^0+\alpha_2x_2^0||
\le \zeta_\nu \le M_{\nu+1}^{-2} \le
 16 (\max \{x_1^0,x_2^0\})^{-2}
.
$$
 
Corollary 1 is proved.

\vskip+0.5cm

Put
$$
A_\nu =\frac{M_\nu^{1/\tau}}{120}.
$$

\vskip+0.5cm

{\bf Corollary 2.}\,\,{\it
Suppose that
\begin{equation}
\zeta_{\nu } \ge A_\nu M_{\nu+1}^{-\frac{\tau}{\tau -1}}
.
\label{stazet}
\end{equation}
 Then there exists an integer point ${\bf x}^0= (x_1^0,x_2^0)$ such that
}

1.\,\, $x_1^0, x_2^0 >0$;

2.\,\, $M_{\nu+1}\le \max \{x_1^0,x_2^0 \} \le  2M_{\nu+1}^{\tau} A_\nu^{-1}$;

3.\,\, $||\alpha_1x_1^0+\alpha_2 x_2^0 ||\le24^\tau M_\nu^{\frac{1-\tau}{\tau}} (\max\{ x_1^0,x_2^0\})^{-\tau}.$

\vskip+0.5cm

Proof.

Apply Lemma 1.
The numbers
 $x_j^0$ from Lemma 1 are positive.
The inequality (\ref{stazet}) leads to the bound
$$
R_\nu \le 2A_\nu^{-1}M_{\nu+1}^{\frac{1}{\tau -1}}  = 2A_\nu^{-1}M_{\nu+1}^\tau
$$
 (as  $\tau^2 =\tau +1$).
This argument in view of  the statement 2 from Lemma 1 together with the Remark to Lemma 1 lead to the statement 2 of Corollary 2.
Moreover, from  (\ref{stazet})  we see that
$$
||\alpha_1x_1^0+\alpha_2x_2^0||\cdot (\max\{x_1^0,x_2^0\})^\tau \le
\zeta_{\nu }R_\nu^\tau
=2^\tau\zeta_{\nu}^{1-\tau} M_{\nu+1}^{-\tau} 
\le 2^\tau A_\nu^{1-\tau}\le
24^\tau M_\nu^{\frac{1-\tau}{\tau}}.
$$
 
Corollary 2 is proved.

\vskip+0.5cm

{\bf Lemma 2.}\,\,{\it Consider consecutive best approximation vectors 
${\bf m}_j,\,\, j = \nu -1,\nu,\nu+1$ such that the inequality  (\ref{neq}) holds.
Suppose that the following two inequalities are valid:
\begin{equation}
\zeta_\nu \le (8M_{\nu-1}M_{\nu+1})^{-1},\,\,\,\,\,
\zeta_{\nu+1} \le (8M_{\nu-1}M_{\nu})^{-1}
.
\label{gather}
\end{equation}
Then there exists an integer point
$(x_1^0,x_2^0)$  such that
}

1.\,\, $x_1^0, x_2^0 >0$;

2.\,\, $\max \{x_1^0,x_2^0 \} \le 20M_{\nu+1} $;

3.\,\, $||\alpha_1x_1^0+\alpha_2 x_2^0 ||< 40 M_{\nu+1}M_\nu^{-1}\zeta_\nu .$

\vskip+0.5cm
Proof.

As
$$
1\neq
\left|
\begin{array}{ccc}
m_{0,\nu-1} &m_{1,\nu -1} & m_{2,\nu-1}\cr
m_{0,\nu} &m_{1,\nu } & m_{2,\nu}\cr
m_{0,\nu+1} &m_{1,\nu +1} & m_{2,\nu+1}
\end{array}
\right|=
\left|
\begin{array}{ccc}
\zeta_{\nu-1} &m_{1,\nu -1} & m_{2,\nu-1}\cr
\zeta_{\nu} &m_{1,\nu } & m_{2,\nu}\cr
\zeta_{\nu+1} &m_{1,\nu +1} & m_{2,\nu+1}
\end{array}
\right|
,
$$
we see that
 $$
1\le
|m_{1,\nu}m_{2,\nu+1}-m_{2,\nu}m_{1,\nu+1}|\zeta_{\nu -1} + 2M_{\nu -1}M_{\nu+1}\zeta_\nu+2M_{\mu-1}M_{\nu}\zeta_{\nu+1}.
$$
We apply   (\ref{gather}) to see that
$$
D_\nu:=
|m_{1,\nu}m_{2,\nu+1}-m_{2,\nu}m_{1,\nu+1}|
\ge (2\zeta_{\nu-1})^{-1}.
$$
From (\ref{mink})  with $\nu$ replaced by $\nu-1$ we have 
\begin{equation}
D_\nu \ge M_\nu^2 /2.
\label{de}
\end{equation}
Consider two-dimensional integer vectors
$$
\xi_\nu = (m_{1,\nu},m_{2,\nu}),\,\,\,\,\,
\xi_{\nu+1} = (m_{1,\nu+1},m_{2\nu+1}),
$$
and the lattice
$$
\Lambda_\nu =\langle \xi_\nu,\xi_{\nu+1}\rangle_{\mathbb{Z}}.$$ 
The fundamental two-dimensional volume of the lattice $\Lambda_\nu$
is equal to 
$ D_\nu $.

Define $\xi_\nu^\perp$
to be the vector of the unit length orthogonal to the vector $\xi_\nu$.
Consider the rectangle $\Omega_\nu^2$ which consists of all points of the form
$$
{\bf x} = \theta_1\xi_\nu +\theta_2 \xi_\nu^\perp,\,\,\,\,\,
|\theta_1 |\le 1,\,\,\, |\theta_2|\le  D_\nu M_{\nu}^{-1}.
$$
Then the measure of $\Omega_\nu^2$  is
$\ge 4D_\nu$.
As $\Omega_\nu^2$  is a convex $0$-symmetric body we may apply the Minkowski convex body theorem.
This theorem ensures that in $\Omega_\nu^2$  there exists a point of  $\Lambda$, independent on $\xi_\nu$.
Hence the rectangle $\Omega_\nu^2$  (as well as any of its translations)
covers a certain fundamental domain with respect to the lattice  $\Lambda$.
Note that the inequality 
  (\ref{de})  leads to the inequality
$ 2D_\nu M_{\nu}^{-1} \ge M_\nu$.
So we see that 
any circle of the radius
 $4D_\nu M_\nu^{-1}$ 
covers a certain fundamental domain with respect to the lattice
 $\Lambda$.
Particulary any circle of the radius 
 $4D_\nu M_\nu^{-1}$  
covers  at least one point of the lattice 
 $\Lambda$.
We take a circle ${\cal C}$ of the radius $4D_\nu M_\nu^{-1}$ centered at the point
$(5D_\nu M_\nu^{-1}, 5D_\nu M_\nu^{-1})$.
Then the point  ${\bf x}^0=(x_1^0,x_2^0)\in {\cal C}\cap \Lambda$
has positive coordinates
 $x_1^0,x_2^0>0$.
Moreover
$$
{\bf x}^0 = \lambda_\nu\xi_\nu+\lambda_{\nu+1}\xi_{\nu+1},
$$
with integer $\lambda_\nu,\lambda_{\nu+1}$.

As 
\begin{equation}
|x_1^0|
=|m_{1,\nu}\lambda_1+m_{1,\nu+1}\lambda_{\nu+1}|\le
10 D_\nu M_\nu^{-1},  
\,\,\,
|x_2^0|
=|m_{2,\nu}\lambda_1+m_{2,\nu+1}\lambda_{\nu+1}|\le
10 D_\nu M_\nu^{-1},  
\label{o}
\end{equation}
and $M_j = \max\{ |m_{1,j}|,|m_{2,j}|\}$,
 we see that
$$
|\lambda_\nu|\le 20 M_{\nu+1} M_\mu^{-1}
,\,\,\,
|\lambda_{\nu+1}|\le 20.
$$
So
$$
||\alpha_1x_1^0+\alpha_2x_2^0||\le
|\lambda_\nu|\zeta_\nu+|\lambda_{\nu+1}|\zeta_{\nu+1}\le
20M_{\nu+1}M_\nu^{-1}\zeta_\nu+20 \zeta_{\nu+1}\le 40M_{\nu+1}M_\nu^{-1}\zeta_\nu.
$$
The statement 3 of Lemma 2 is proved.

As
 $D_\nu \le 2M_\nu M_{\nu+1}$ we deduce from   (\ref{o}) 
that
$$
\max\{ x_1^0, x_2^0\} \le 20 M_{\nu+1}.
$$
So the statement 2 of Lemma 2 is verified and Lemma 2 is proved.

\vskip+0.5cm
{\bf Remark.}\,\, The contitions (\ref{gather})
are technical. Unfortunately we cannot avoid them as the inequalites (\ref{mink}) with $j=\nu,\nu+1$ are not sufficient for the proof.
 
\vskip+0.5cm
{\bf Corollary 3.}\,\,{\it
 Let the inequality  (\ref{neq}) be valid for $\nu $  large enough.
 Suppose that
\begin{equation}
\zeta_\nu <A_\nu M_{\nu+1}^{-\frac{\tau}{\tau -1}}.
\label{n1}
\end{equation}
In addition suppose that the second inequality from  (\ref{gather}) is also  satisfied.
Then
there exists an integer point
 $(x_1^0,x_2^0)$ such that
}

1.\,\, $x_1^0, x_2^0 >0$;

2.\,\, $\max \{x_1^0,x_2^0 \} \le 20M_{\nu+1} $;

3.\,\, $||\alpha_1x_1^0+\alpha_2 x_2^0 ||< 24^\tau M_\nu^{\frac{1-\tau}{\tau}} (\max\{x_1^0,x_2^0\})^{-\tau} .$

\vskip+0.5cm
Proof.

For $\nu$ large enough   the first inequality   from (\ref{gather}) follows from (\ref{n1}).
Now Corollary 3 immediately follows from Lemma 2.

 \section{Proof of Theorem 3}

Suppose that
 $\zeta_{\nu+1}\ge (8M_{\nu+2}^2)^{-1}$. Then
we apply Corollary 1 (with $\nu$ replaced by $\nu+1$).
The statement {\bf (i)} of Theorem 3 follows.

Suppose that
  $\zeta_{\nu+1} < (8M_{\nu+2}^2)^{-1}$.  As
 $M_{\nu+2}\ge M_\nu\ge M_{\nu -1}$  and $\zeta_{\nu+1} \le M_{\nu+2}^{-2}$ we see that
the second inequality from  (\ref{gather}) is satisfied. 
In the case  $\zeta_\nu \ge A_\nu M_{\nu+1}^{\frac{\tau}{\tau -1}}$
 we apply Corollary 2. 
In the case  $\zeta_\nu \le A_\nu M_{\nu+1}^{\frac{\tau}{\tau -1}}$
 we apply Corollary 3. 
 So we establish statements  1, 3 from {\bf (ii)}.
The statement 2 from {\bf (ii)} also follows from the Corollaries 2,3 as  $M_{\nu+1}^{\tau -1} = M_{\nu+1}^{1/\tau}\ge M_\nu^{1/\tau}$.

Theoreem 3 is proved.

author:

\noindent Nikolay  G. {\sc Moshchevitin} \\
Moscow Lomonosov State University \\
Vorobiovy Gory, GSP--1 \\
119991 Moscow, RUSSIA \\
\emph{E-mail}: {\fontfamily{cmtt}\selectfont moshchevitin@mech.math.msu.su,
moshchevitin@rambler.ru}

\begin{thebibliography}{100}
\bibitem{SCH}
W.M. Schmidt,\,{\it
Two questions in Diophantine approximations},\,
Monatshefte f\"ur Mathematik\, {\bf 82}, 237 - 245 (1976). 

\bibitem{SCH1}
  W. M. Schmidt,\, {\it
 Open problems in Diophantine approximations}, in   "Approximations Diophantiennes et nombres transcendants' Luminy, 1982, Progress in
 Mathematics,
 Birkh\"auser, p.271 - 289 (1983).




\bibitem{T1}
P. Thurnheer,\,
{\it
Zur diophantischen Approximation von zwei reellen Zahlen},
       Acta Arithmetica  {\bf 44}, 201-206 (1984).



\bibitem{T2}
P. Thurnheer,\,
{\it On Dirichlet's theorem concerning Diophantine approximation},
    Acta Arithmetica {\bf 54}, 241-250 (1990).

\bibitem{BK}
Y. Bugeaud, S. Kristensen,\,{\it
Dipohantine exponents for mildly restricted approximation}, Preprint avaliable at: arXiv
\bibitem{DS}
H. Davenport, W.M. Schmidt,\,{\it
Approximation to real numbers by quadratic irrationals},\,  Acta Arithmetica {\bf  13},  169 - 176 (1967).

\bibitem{M}
 N.G.Moshchevitin,\, {\it Best Diophantine approximations: the phenomenon of
    degenerate dimension.} London Mathematical Society, Lecture Note Series {\bf 338},
    158--182 (2007). 


\end{thebibliography}
\end{document}